\documentclass[11pt,reqno]{amsart}
\usepackage{graphicx}
\usepackage{amssymb,latexsym,epsfig,mathrsfs,graphpap,graphics,amsmath,amssymb}
\usepackage{amstext}
\usepackage{verbatim}
\usepackage{color}

\newtheorem{conj.}[thm]{Conjecture}

\theoremstyle{definition}

\theoremstyle{remark}

\numberwithin{equation}{section}

\begin{document}

\begin{flushleft}
  {\bf\Large { Special Affine  Wavelet  Transforms and the Corresponding\\[1.6mm] Poisson Summation Formula}}
\end{flushleft}

\parindent=0mm \vspace{.4in}

  {\bf{Firdous A. Shah$^{\star}$, Azhar Y. Tantary$^{\star}$ and Aajaz A. Teali$^{\star}$ }}

 \parindent=0mm \vspace{.1in}
{\small \it $^{\star}$Department of  Mathematics,  University of Kashmir, South Campus, Anantnag-192101, Jammu and Kashmir, India. E-mail: $\text{fashah@uok.edu.in}$; $\text{aytku92@gmail.com}$;$\text{aajaz.math@gmail.com}$}

\parindent=0mm \vspace{.2in}
{\small{\bf Abstract.} The special affine Fourier transform (SAFT) is a promising tool for analyzing non-stationary signals with more degrees of freedom.  However, the SAFT fails in obtaining the local features of non-transient signals due to its global kernel and thereby make SAFT incompetent in situations demanding joint information of time and frequency. To circumvent this limitation, we propose a highly flexible time-frequency transform namely, the special affine wavelet transform (SAWT) and investigate the associated constant $Q$-property in the joint time-frequency domain. The basic properties of the proposed transform such as Rayleigh's theorem, inversion formula  and characterization of the range are discussed using the machinery of special affine Fourier transforms and operator theory. Besides this, the discrete counterpart of SAWT is also discussed and the corresponding reconstruction formula is obtained. Moreover, we also drive a direct relationship between the well known special affine Wigner distribution and the proposed transform. This is followed by introducing a new kind of wave packet transform associated with the special affine Fourier transform. Towards the end, an analogue of the Poisson summation formula for the proposed special affine wavelet transform is derived.

\parindent=0mm \vspace{.1in}

\parindent=0mm \vspace{.1in}
}

\section{Introduction}

During the culminating years of last century, a reasonably flexible integral transform associated with a general inhomogeneous lossless linear mapping in phase-space was introduced, namely, the special affine Fourier transform  \cite{AS,AS1,Cai}. The SAFT is a six-parameter class of linear integral transform which encompasses a number of well known unitary transforms including the classical Fourier transform, fractional Fourier transform, Fresnel transform, Laplace transform, Gauss-Weierstrass transform, Bargmann transform and the linear canonical transform \cite{DS,DS1,Nv,Al,JA,HKOS}.  For a uni-modular matrix $M=\left(A,B,C,D, p, q\right)$, the special affine Fourier transform of any signal $f$ with respect to the matrix $M$ is defined by
\begin{align*}\label{1.1}
\mathcal O^M\big[f\big](a)=\int_{-\infty}^{\infty}f(t)\,\mathcal{K}_M(t,a)\,dt,\tag{1.1}
\end{align*}
where the kernel $\mathcal{K}_M(t,a)$  is given by
\begin{align*}
\mathcal{K}_M(t,a)=\dfrac{1}{\sqrt {2\pi iB}}\exp\left\{\dfrac{i}{2B}\Big(At^2+2t(p-a)-2a(Dp-Dq)+D(a^2+p^2)\Big) \right\}.\tag{1.2}
\end{align*}

It is pertinent to note that  the SAFT is just a chirp multiplication operation if $B=0$, therefore, for brevity this case is omitted in theoretical studies. We also note that the  phase-space transform (1.1) is lossless if and only if $|\det M|=1$, that is; $AD-BC=1$ and for this reason, SAFT is also known as the inhomogeneous canonical transform \cite{AS1}. Due to the extra degrees of freedom, SAFT  has attained a respectable status within a short span and is being broadly employed across  several disciplines of science and engineering including  signal and image processing, optical and radar systems, electrical and communication systems, pattern
recognition, sampling theory, shift-invariant theory and quantum mechanics \cite{Qzy,XQ,ZWZ,Z,Xhch,ULT,BZ}. Despite the versatile applicability, the SAFT surfers from an undeniable drawback  due to its global kernel involved in (1.1) and thereby is incompetent in situations demanding a joint analysis of time and spectral characteristics of a signal.  To address this issue, an immediate concern is to revamp the existing special affine Fourier transform by adjoining certain excellent localization features. The most appropriate candidates are the wavelet functions, which inherits nice localization properties along with additional characteristic features such as orthogonality, vanishing moments and self adjustability. By intertwining the ideas of SAFT and wavelet transforms, we introduce a hybrid integral transform namely, special affine wavelet transform which is capable of providing a joint time and frequency localization of non-stationary signals with more degrees of freedom.

\parindent=8mm \vspace{.1in}
The main contributions of the article are mentioned below:

\begin{itemize}

 \item To introduce a hybrid integral transform and discuss the associated constant $Q$-property in the joint time-frequency domain.

  \item To study some fundamental properties of the proposed transform including Rayleigh's theorem, inversion formula and characterization of range.
 \item To introduce the discrete counterpart of the proposed transform and obtain the associated inversion formula.

  \item To obtain a direct relationship between the special affine Wigner distribution and the proposed transform.

  \item To introduce a new wave packet transform associated with the special affine Fourier transform.

  \item To obtain an  analogue of the well known Poisson summation formula for the proposed special affine wavelet transform.

\end{itemize}

\parindent=0mm \vspace{.1in}
The rest of the article is organized as follows. The new hybrid transform and some of its basic properties are studied in  Sections 2 and 3. Section 4 is devoted to introducing the discrete counterpart of the special affine wavelet transform. Relationship between the special affine Wigner distribution and the proposed transform is presented in Section 5.  In Section 6, we study the wave packet transform associated with the special affine Fourier transform. The final section is devoted to establishing the Poisson summation formula for the special affine wavelet transform.

\pagestyle{myheadings}

\section{ Continuous Special Affine Wavelet Transform And the Associated Constant $Q$-Property}

In this Section, our aim is to formally introduce the continuous special affine wavelet transform by intertwining the advantages of the special affine Fourier  and wavelet transforms.

\parindent=8mm \vspace{.1in}
Wavelets are arguably the most appropriate mathematical entities with affine-like structure of well-localized waveforms at various scales and locations; that is, they are generated by dilation and translation of a single  mother function  $\psi\in L^2(\mathbb R)$. Mathematically, the wavelet family $\psi_{b,a}(t)$ is defined as
\begin{align*}
{\psi}_{a,b}(t)= \dfrac{1}{\sqrt{|a|}}\, \psi\left(\dfrac{t-b}{a}\right),\quad a,b\in\mathbb{R},\quad a\neq{0},\tag{2.1}
\end{align*}
where $a$ is a scaling parameter which measures the degree of compression or scale, and $b$ is a translation parameter which determines the location of the wavelet \cite{DS}. With major modifications of the family (2.1), we define a new family of functions $\psi^M_{a,b}(t)$ by the combined action of the unitary operators $\mathcal D_a, \mathcal T_b$ and the uni-modular matrix $M=(A,B,C,D,p,q)$ as
\begin{align*}\label{2.1}
\psi^M_{a,b}(t)= \dfrac{1}{\sqrt {a}} \psi\left( \dfrac{t-b}{a} \right)\, \mathcal{K}^M(t,a),\tag{2.2}
\end{align*}
where
\begin{align*}\label{2.2}
\mathcal{K}^M(t,a)= \dfrac{1}{i\sqrt{2\pi iB}}\exp \left\{ \dfrac{-i}{2B}\Big( At^2+2t(p-a)-2a(Dp-Bq)+D(a^2+p^2)\Big) \right\}.\tag{2.3}
\end{align*}

\parindent=0mm \vspace{.0in}
Having formulated a family of analyzing functions, we are now ready to introduce the definition of the continuous special affine wavelet transform in $L^2(\mathbb R)$.

\parindent=0mm \vspace{.1in}
{\bf Definition 2.1.} For any finite energy signal $f\in L^2(\mathbb R)$, the continuous special affine wavelet transform of $f$ with respect to the wavelet  $\psi\in L^2(\mathbb R)$ is defined by
\begin{align*}\label{2.3}
{\mathcal{W}}_\psi^M \big[f\big](a,b)=\Big\langle f,\, \psi^M_{a,b} \Big\rangle
=\int_{-\infty}^{\infty} f(t)\,\overline{\psi^{M}_{a,b}(t)}\,dt.\tag{2.4}
\end{align*}
where $\psi^M_{a,b}(t)$ is given by (2.1).

\parindent=8mm \vspace{.1in}
It is worth noting that the proposed transform (2.4) boils down to some existing integral transforms as well as gives birth to some new time-frequency transforms as mentioned below:

\parindent=0mm \vspace{.1in}
(i). For $M=\left(A,B,C,D,0,0\right)$, the continuous special affine wavelet transform (2.4) reduces to the  linear canonical wavelet transform given by \cite{WL,WWWR}
\begin{align*}\label{2.5}
\left[{\mathcal{W}}_\psi^M f\right](a,b)=\dfrac{1}{\sqrt {2a\pi iB}}\int_{-\infty}^{\infty}f(t)\,\overline{\psi\left( \dfrac{t-b}{a} \right)}\exp\left\{\dfrac{i\big(At^2-2ta+Da^2\big)}{2B} \right\}dt.\tag{2.5}
\end{align*}

\parindent=0mm \vspace{.0in}
(ii). For the matrix $M=\left(\cos\theta,\sin\theta,-\sin\theta,\cos\theta,0,0\right)$, the continuous special affine wavelet transform (2.4) boils down to the  fractional wavelet transform defined  by \cite{DZW}
      \begin{align*}\label{2.6}
     \left[{\mathcal{W}}_\psi^M  f\right](a,b)=\dfrac{1}{\sqrt {2a\pi i\sin \theta}}\int_{-\infty}^{\infty}f(t)\,\overline{\psi\left( \dfrac{t-b}{a} \right)}
        \exp\left\{\dfrac{i\big(t^2+a^2\big)\cot\theta }{2}-ita\, \csc\theta\right\}dt.\tag{2.6}
     \end{align*}

\parindent=0mm \vspace{.0in}
(iii). For the matrix $M=\left(\cos\theta,\sin\theta,-\sin\theta,\cos\theta,p,q\right)$, we can obtain a new transform called the  special affine fractional wavelet  transform defined  by
      \begin{align*}\label{2.7}
     \left[{\mathcal{W}}_\psi^M  f\right](a,b)&=\dfrac{1}{\sqrt {2a\pi i\sin \theta}}\int_{-\infty}^{\infty}f(t)\,\overline{\psi\left( \dfrac{t-b}{a} \right)} \\
        &\quad\times \exp\left\{\dfrac{i\big(t^2+a^2+p^2-2ap\big)\cot\theta }{2}+it(p-a)\, \csc\theta+iaq\right\}\,dt.\tag{2.7}
     \end{align*}

\parindent=0mm \vspace{.0in}
(iv). For the matrix $M=\big(1,iB,0,1,p,q\big)$,    we can obtain a new transform namely the special affine Gauss-Weierstrass-wavelet transform given by
  \begin{align*}\label{2.8}
   \left[{\mathcal{W}}_\psi^M f\right](a,b)=\dfrac{1}{\sqrt {2a\pi (-B)}}\int_{-\infty}^{\infty}f(t)\overline{\psi\left( \dfrac{t-b}{a} \right)} \exp\left\{\dfrac{ (t-a)^2+2tp+p^2-2a(p-iBq)}{2B}\right\}dt.\tag{2.8}
     \end{align*}

\parindent=0mm \vspace{.0in}
(v). For the matrix $M=(1,B,0,1,p,q)$,    we can obtain a new transform namely the special affine Fresnel-wavelet transform given by
  \begin{align*}\label{2.9}
    \left[{\mathcal{W}}_\psi^M  f\right](a,b)&=\dfrac{1}{\sqrt {2a\pi i B}}\int_{-\infty}^{\infty}f(t)\,\overline{\psi\left( \dfrac{t-b}{a} \right)} \exp\left\{\dfrac{i\big( (t-a)^2+2tp+p^2-2a(p-Bq)\big)}{2B}\right\}\,dt.\tag{2.9}
     \end{align*}

\parindent=0mm \vspace{.0in}
(vi).  By changing the usual matrix $M=\big(A,B,C,D,p,q\big)$ to  $M=(1,-iB,0,1,p,q),$ we can obtain special affine Bilateral Laplace-wavelet transform defined by
    \begin{align*}\label{2.10}
   \left[{\mathcal{W}}_\psi^M  f\right](a,b)=\dfrac{1}{\sqrt {2a\pi B}}\int_{-\infty}^{\infty}f(t)\overline{\psi\left( \dfrac{t-b}{a} \right)}\exp\left\{-\dfrac{ (t-a)^2+2tp+p^2-2a(p+iBq)}{2B}\right\}dt.\tag{2.10}
     \end{align*}

\parindent=0mm \vspace{.1in}
For a lucid demonstration of the continuous special affine wavelet transform (2.4), we shall present some illustrative examples.

\parindent=0mm \vspace{.1in}

{\bf Example 2.2} (i). Consider the constant function $f(t)=K$ and the Morlet function $\psi(t)=\exp\big\{i\alpha t-\frac{t^2}{2} \big\}$.
Then, the translated and scaled versions of $\psi(t)$ are given by
\begin{align*}
\psi\left(\frac{t-b}{a}\right)= \exp\left\{i\alpha \left(\frac{t-b}{a}\right)-\frac{t^2+b^2-2tb}{2a^2}\right\}.
\end{align*}
Consequently, the family of analyzing functions $\psi^M_{a,b}(t)$ is obtained as
\begin{align*}
\psi^M_{a,b}(t)= \dfrac{1}{\sqrt {2\pi iBa}}\psi\left(\dfrac{t-b}{a}\right)\exp\left\{ \dfrac{-i}{2B}\Big( At^2+2t(p-a)-2a(Dp-Bq)+D(a^2+p^2)\Big) \right\}.
\end{align*}
To compute the continuous special affine wavelet transform of $f(t)$ with respect to the uni-modular matrix $M=\left(A,B,C,D, p, q\right)$ and the window function $\psi(t)$, we proceed as
\begin{align*}
&\Big[{\mathcal{W}}_\psi^M  f\Big](a,b)\qquad\\
&=\dfrac{K}{\sqrt {2\pi iBa}}\exp \left\{ \frac{i}{2B}\Big(-2a(Dp-Bq)+D(a^2+p^2)\Big) \right\}\\
&\qquad\qquad\times\int_{-\infty}^{\infty} \exp\left\{\frac{iAt^2+2t(p-a)}{2B}\right\} \overline{ \exp\left\{i\alpha \left(\frac{t-b}{a}\right)-\frac{t^2+b^2-2tb}{2a^2}\right\}}\,dt\\
&=\dfrac{K}{\sqrt {2\pi iBa}}\exp \left\{ \frac{i}{2B}\Big(-2a(Dp-Bq)+D(a^2+p^2)\Big)+\frac{ib\alpha}{a}-\frac{b^2}{2a^2} \right\}\\
&\qquad\qquad\qquad\times\int_{-\infty}^{\infty} \exp\left\{-t^2\Big(\frac{1}{2a^2}-\frac{iA}{2B} \Big)+t\Big(\frac{p-a}{B}-\frac{i\alpha}{a}+\frac{b}{a^2} \Big)\right\}\,dt\\
&=\dfrac{K}{\sqrt {2\pi iBa}}\exp \left\{ \frac{i}{2B}\Big(-2a(Dp-Bq)+D(a^2+p^2)\Big)+\frac{ib\alpha}{a}-\frac{b^2}{2a^2} \right\}\\
&\qquad\qquad\qquad\times \sqrt{\frac{\pi 2Ba^2}{B-iAa^2}} \exp\left\{\Big(\dfrac{(p-a)a^2-i\alpha aB+bB}{a^2B}\Big)^2\times \frac{2a^2B}{4(b-iAa^2)}\right\}\\
&= \frac {K \sqrt{a}}{\sqrt{iB+a^2A}} \exp \left\{ \frac{i}{2B}\Big(-2a(Dp-Bq)+D(a^2+p^2)\Big)+\frac{ib\alpha}{a}-\frac{b^2}{2a^2} \right\}\\
&\qquad\qquad\qquad\qquad\qquad\qquad\qquad\qquad\qquad\times\exp\left\{\frac{\Big((p-a)a^2-i\alpha aB+bB \Big)^2}{2a^2B(B-ia^2A)}\right\}.
\end{align*}

(ii). For the  Haar function
\begin{align*}
\psi(t)=\left\{\begin{array}{cc}
1, & 0\le t < 1/2 \\
-1, & 1/2 \le t\le 1 ,
\end{array}\right.\tag{2.11}
\end{align*}
we first construct the family of analyzing functions $\psi^M_{a,b}(t)$ in the following way:
\begin{align*}
&\psi^M_{a,b}(t)\\
&=\left\{\begin{array}{cc}
\dfrac{1}{\sqrt {2\pi iBa}}\exp\left\{ \dfrac{-i}{2B}\Big( At^2+2t(p-a)-2a(Dp-Bq)+D(a^2+p^2)\Big) \right\}, & b\le t < \dfrac{a}{2}+b \\
-\dfrac{1}{\sqrt {2\pi iBa}}\exp\left\{ \dfrac{-i}{2B}\Big( At^2+2t(p-a)-2a(Dp-Bq)+D(a^2+p^2)\Big) \right\}, & \dfrac{a}{2}+b \le t\le a+b .
\end{array}\right.
\end{align*}
Then, we compute the continuous special affine wavelet transform of the generalized Gaussian function $f(t)=e^{-iAt^2/2B},AB>0$  in the following manner:

\begin{align*}
&\Big[{\mathcal{W}}_\psi^M  f\Big](a,b)\\
&=\dfrac{1}{\sqrt {2\pi iBa}}\exp \left\{ \dfrac{i}{2B}\Big(-2a(Dp-Bq)+D(a^2+p^2)\Big) \right\} \\
&\qquad\qquad \times \Bigg\{\int_{b}^{\frac{a}{2}+b} \exp \left\{ \dfrac{i}{2B}\Big(2t(p-a)\Big) \right\}\,dt  - \int_{\frac{a}{2}+b}^{a+b} \exp \left\{ \dfrac{i}{2B}\Big(2t(p-a)\Big) \right\}\,dt \Bigg\}\\
&=\dfrac{1}{\sqrt {2\pi iBa}}\exp \left\{ \dfrac{i}{2B}\Big(-2a(Dp-Bq)+D(a^2+p^2)\Big) \right\}\, \dfrac{B}{i(p-a)} \Bigg[\exp \left\{ \dfrac{i(p-a)}{B}\left(\frac{a}{2}+b\right) \right\}\\
&\qquad\qquad\qquad - \exp \left\{ \dfrac{i(p-a)b}{B} \right\}- \exp \left\{ \dfrac{i(p-a)(a+b)}{B} \right\}+\exp \left\{ \dfrac{i(p-a)}{B}\left(\frac{a}{2}+b\right) \right\} \Bigg]\\
&=\dfrac{\sqrt{iB}}{\sqrt {2\pi a(a-p)}}\exp \left\{ \dfrac{i}{2B}\Big(-2a(Dp-Bq)+D(a^2+p^2)\Big) \right\}\, \\
&\qquad\qquad\qquad\qquad\qquad\qquad\times \exp \left\{ \dfrac{ib(p-a)}{B} \right\} \Bigg[ 2 \exp \left\{ \dfrac{ia(p-a)}{2B}\right\}-1- \exp \left\{ \dfrac{i(p-a)a}{B}\right\}\Bigg].
\end{align*}

\parindent=8mm \vspace{.0in}
In the following, we shall derive a fundamental relationship between the continuous special affine wavelet transform \eqref{2.3} and the special affine Fourier transform \eqref{1.1}. With the aid of this formula, we shall study the $Q$-property of the proposed  transform and obtain the associated joint time-frequency resolution.

\parindent=0mm \vspace{.1in}

{\bf Proposition 2.3.}  {\it Let ${\mathcal{W}}_\psi^M  \big[f\big](a,b)$ and ${\mathcal{O}}^{M}\big[f\big](a)$  be the continuous  special affine wavelet transform  and the special affine Fourier transform of any finite energy signal $f\in L^2(\mathbb R)$, respectively. Then, we have }
\begin{align*}\label{2.6}
&{\mathcal{W}}^{M}_{\psi}\big[f\big](a,b)\\
&=\sqrt{a} \int_{-\infty}^{\infty} \exp\left\{\dfrac{-i}{2B}\Big( 2a(\xi -1)(Dp-Bq)+Da^2(1-\xi^2)-2B(p-a)\Big)\right\}\,{\mathcal{O}}^{M}\big[f\big](\xi)\\
&\times \overline{{\mathcal{O}}^{M}}\left[\psi(y) \exp\left\{\dfrac{i}{2B}\Big((2a^2-1)Ay^2-2yp+2ay(2Ab-2p+\xi-a)\Big)\right\} \right](a\xi) \overline{\mathcal{K}^M(b,\xi)}\,d\xi \tag{2.12}
\end{align*}
where $\mathcal{K}^{M}(b,\xi)$ is given by \eqref{2.2}.

\parindent=0mm \vspace{.1in}
{\it Proof.} Applying the definition of special affine Fourier transform, we obtain
\begin{align*}
&\mathcal{O}^M \Big[ \psi^M_{a,b}(t)\Big](\xi)\\
&= \int_{-\infty}^{\infty} \psi^M_{a,b}(t)\, \mathcal{K}^M(t,\xi)dt\\
&=\int_{-\infty}^{\infty} \dfrac{1}{\sqrt{a}} \psi \left( \dfrac{t-b}{a}\right) \dfrac{1}{\sqrt{2\pi iB}} \exp\left\{\dfrac{i}{2B}\Big( At^2+2t(p-a)-2a(Dp-Bq)+D(a^2+p)\Big) \right\}\\
&~\times\dfrac{1}{\sqrt{2\pi iB}} \exp\left\{\dfrac{i}{2B}\Big( At^2+2t(p-\xi)-2a(Dp-Bq)+D(\xi^2+p)\Big) \right\} dt\\
&=\int_{-\infty}^{\infty} \dfrac{1}{\sqrt{a}} \psi(y) \dfrac{1}{2\pi iB} \exp\left\{\dfrac{i}{2B}\Big( A(ay+b)^2+2(ay+b)(p-a)-2a(Dp-Bq)+D(a^2+p)\Big) \right\}\\
&~\times\exp\left\{\dfrac{i}{2B}\Big( A(ay+b)^2+2(ay+b)(p-\xi)-2a(Dp-Bq)+D(\xi^2+p)\Big) \right\} \, a \, dy\\
\end{align*}
\begin{align*}
&=\int_{-\infty}^{\infty} \psi(y) \dfrac{1}{2\pi iB\sqrt{a}} \exp\left\{\dfrac{i}{2B}\Big( A(a^2y^2+b^2+2ay)+2(ay+b)(p-a)-2a(Dp-Bq)+D(a^2+p)\Big) \right\}\\
&~\times\exp\left\{\dfrac{i}{2B}\Big( A(a^2y^2+b^2+2ay)+2(ay+b)(p-\xi)-2a(Dp-Bq)+D(\xi^2+p)\Big) \right\} a\, dy\\
&=\dfrac{\sqrt{a}}{\sqrt{2\pi iB}} \int_{-\infty}^{\infty}  \psi(y) \dfrac{1}{\sqrt{2\pi iB}} \exp\left\{\dfrac{i}{2B}\Big( Ab^2+2b(p-\xi)-2\xi(Dp-Bq)+D(\xi^2+p)\Big) \right\}\\
&~\times\exp\left\{\dfrac{i}{2B}\Big( 2Aa^2y^2+4Aaby+2ay(p-a)+2ay(p-a)+Ab^2+2b(p-a)\Big) \right\}\\
&~\times \exp\left\{\dfrac{i}{2B}\Big(-2a(Dp-Bq)+D(\xi^2+p)\Big) \right\} \, dy\\
&=\dfrac{\sqrt{a}}{\sqrt{2\pi iB}} \int_{-\infty}^{\infty}  \psi(y) \, \mathcal{K}^M(b,\xi)\exp\left\{\dfrac{i}{2B}\Big( 2Aa^2y^2-Ay^2+Ab^2\Big) \right\}\\
&~ \times \exp\left\{\dfrac{i}{2B}\Big( Ay^2+2y\big(p-(a\xi)\big)-2(a\xi)(Dp-Bq)+D(a^2\xi^2+p)\Big) \right\}\\
&~\times\exp\left\{\dfrac{i}{2B}\Big(4Aaby+2ay(p-a)+4ayp-2a^2y+2bp-2ba-2a(Dp-Bq)+Da^2\Big) \right\}\\
&~\times\exp\left\{\dfrac{i}{2B}\Big(-2yp+2ya\xi+2a\xi(Dp-Bq)-D(a\xi)^2\Big) \right\} \, dy\\
&=\sqrt{a}\, \mathcal{K}^M(b,\xi)\, \exp \left\{ \dfrac{i}{2B}\Big(2bp-2ab-2a(Dp-Bq)+Da^2+2a\xi(Dp-Bq)-Da^2\xi^2+Ab^2\Big) \right\}\\
&~\times \dfrac{1}{\sqrt{2\pi iB}} \int_{-\infty}^{\infty} \Big[ \psi(y)\exp\left\{\dfrac{i}{2B}\Big(2Aa^2y^2-Ay^2+4Aaby+4ayp-2yp+2ay\xi-2a^2y\Big) \right\}\Big]\\
&~\times \exp\left\{\dfrac{i}{2B}\Big( Ay^2+2y[p-(a\xi)]-2(a\xi)(Dp-Bq)+D(a^2\xi^2+p)\Big) \right\}dy\\
&=\sqrt{a}\, \exp\left\{\dfrac{i}{2B}\Big(2bp-2ab-2a(Dp-Bq)+Da^2+Ab^2+2a\xi(Dp-Bq)-Da^2\xi^2\Big) \right\}\\
&~ \times \mathcal{O}^M \Big[ \psi(y)\exp\left\{\dfrac{i}{2B}\Big(2Aa^2y^2-Ay^2+4Aaby+4ayp-2yp+2ay\xi-2a^2y\Big) \right\}\Big](a\xi)\, \mathcal{K}^M(b,\xi)\\
&=\sqrt{a}\, \exp\left\{\dfrac{i}{2B}\Big(2a(\xi-1)(Dp-Bq)+Da^2(1-\xi^2)-2b(p-a)+Da^2(1-\xi^2)+Ab^2\Big) \right\}\\
&~ \times \mathcal{O}^M \Big[ \psi(y)\exp\left\{\dfrac{i}{2B}\Big((2a^2-1)Ay^2-2yp+2ay(2Ab-2p+\xi-a)\Big) \right\}\Big](a\xi)\, \mathcal{K}^M(b,\xi). \tag{2.13}
\end{align*}
By virtue of  (2.13), the special affine wavelet transform \eqref{2.3} can be expressed as
\begin{align*}
&{\mathcal{W}}^{M}_\psi \big[f\big](a,b)\\
&~=\Big\langle \mathcal{O}^M \big[f\big],~\mathcal{O}^M\left[ \psi^M_{a,b}\right] \Big\rangle\\
&~=\int_{-\infty}^{\infty} \mathcal{O}^M \big[f\big](\xi)\, \overline{\mathcal{O}^M}\left[ \psi^M_{a,b}\right](\xi)\,d\xi\\
&~=\sqrt{a}\,\int_{-\infty}^{\infty} \overline{ \mathcal{K}^M(b,\xi)}\, \mathcal{O}^M \big[f\big](\xi)\, \exp\left\{\dfrac{-i}{2B}\Big(2a(\xi-1)(Dp-Bq)+Da^2(1-\xi^2)-2b(p-a)\Big) \right\}\\
&\qquad\times \overline {\mathcal{O}^M} \left[ \psi(y)\exp\left\{\dfrac{i}{2B}\Big((2a^2-1)Ay^2-2yp+2ay(2Ab-2p+\xi-a)\Big) \right\}\right](a\xi)\\
&\qquad \times \exp\left\{\dfrac{-i}{2B}\Big(Da^2(1-\xi^2)+Ab^2\Big) \right\}\,d\xi.
\end{align*}

\parindent=0mm \vspace{.0in}
This completes the proof of the Proposition 2.3. \quad\fbox

\parindent=8mm \vspace{.1in}
As a consequence of Proposition 2.3, we conclude that if the analyzing functions $\psi^M_{a,b}(t)$ are supported in the time-domain or the special affine Fourier domain, then the proposed transform $\mathcal{W}^M_\psi \big[f \big](a,b)$ is accordingly supported in the respective domains. This implies that the special affine wavelet transform is capable of providing the simultaneous information of the time and the special affine frequency in the time-frequency domain. To be more specific, suppose that $\psi(t)$ is the window with centre $E_{\psi}$ and radius $\Delta_{\psi}$ in the time domain. Then, the centre and radii of the time-domain window function $\psi_{a,b}^{M}(t)$  of the proposed transform (2.4) is given by
\begin{align*}\label{2.6}
E\Big[\psi^M_{a,b}(t)\Big]&=\dfrac{\displaystyle\int_{-\infty}^{\infty}t\Big|\psi^M_{a,b}(t)\Big|^2 dt}{\displaystyle\int_{-\infty}^{\infty}\Big|\psi^M_{a,b}(t)\Big|^2 dt}=\dfrac{\displaystyle\int_{-\infty}^{\infty}t\Big|\psi_{a,b}(t)\Big|^2 dt}{\displaystyle\int_{-\infty}^{\infty}\Big|\psi_{a,b}(t)\Big|^2dt}=E\Big[\psi_{a,b}(t)\Big]=b+aE_{\psi}\tag{2.14}
\end{align*}
and
\begin{align*}\label{2.7}
\Delta\Big[\psi^M_{a,b}(t)\Big]&=\left(\dfrac{\displaystyle\int_{-\infty}^{\infty}\Big(t-(b+aE_{\psi})\Big)\Big|\psi^M_{a,b}(t)\Big|^2 dt}{\displaystyle\int_{-\infty}^{\infty}\left|\psi^M_{a,b}(t)\right|^2 dt}\right)^{1/2}=\left(\dfrac{\displaystyle\int_{-\infty}^{\infty}\Big(t-b-aE_{\psi}\Big)\Big|\psi_{a,b}(t)\Big|^2 dt}{\displaystyle\int_{-\infty}^{\infty}\Big|\psi_{a,b}(t)\Big|^2dt}\right)^{1/2}\\
&=\Delta\Big[\psi_{a,b}(t)\Big]= a\Delta_{\psi},\tag{2.15}
\end{align*}
respectively. Let $H(\xi)$ be the window function in the special affine Fourier transform domain given by
\begin{align*}
H(\xi)=\mathcal{O}^M\left[\exp\left\{\dfrac{i}{2B}\Big((2a^2-1)Ay^2-2yp+2ay(2Ab-2p+\xi-a)\psi(y)\Big) \right\}\right]\left(\xi\right).
\end{align*}

\parindent=0mm \vspace{.0in}
Then, we can derive the center and radius of the special affine Fourier domain window function
\begin{align*}
H(a\xi)=\mathcal{O}^M\left[\exp\left\{\dfrac{i}{2B}\Big((2a^2-1)Ay^2-2yp+2ay(2Ab-2p+\xi-a)\psi(y)\Big) \right\}\right]\left(a\xi\right)
\end{align*}
appearing in (2.12) as
\begin{align*}
E\big[H\left(a\xi\right)\big]= \dfrac{\displaystyle\int_{-\infty}^{\infty}(a\xi)\big|H(a\xi)\big|^2 d\xi}{\displaystyle\int_{-\infty}^{\infty}\big|H(a\xi)\big|^2 d\xi} = a \, E_H,\tag{2.16}
\end{align*}
and
\begin{align*}
\Delta\big[H\left(a\xi \right)\big] =a \,\Delta_H. \tag{2.17}
\end{align*}
Thus, the $Q$-factor of the  proposed transform (2.4) is given by
\begin{align*}\label{2.10}
Q=\dfrac{\text{width of the window function}}{\text{centre of the window function}}=\dfrac{\Delta\big[H\left(a\xi \right)\big]}{E\big[H\left(a\xi \right)\big]}=\dfrac{\Delta_H}{E_H}=\text{constant},\tag{2.18}
\end{align*}
which is independent of the uni-modular matrix  $M=(A,B,C,D,p,q)$ and the scaling parameter $a$.  Therefore, the localized time and frequency characteristics of the proposed transform (2.4) are given in the  time and frequency windows
\begin{align*}\label{2.13}
\Big[b+aE_{\psi}-a\Delta_{\psi},~b+aE_{\psi}+a\Delta_{\psi}\Big] \quad \text{and}\quad \Big[a\,E_H-a\,\Delta_H,~a\,E_H+a\,\Delta_H\Big],\tag{2.19}
\end{align*}
respectively. Hence, the joint resolution of the continuous special affine wavelet transform (2.4) in the time-frequency domain is described by a flexible window $\psi$ having a total spread  $4\Delta_{\psi}\Delta_H$ and is given by
\begin{align*}\label{2.13}
\Big[b+aE_{\psi}-a\Delta_{\psi},~b+aE_{\psi}+a\Delta_{\psi}\Big]
\times \Big[a\,E_H-a\,\Delta_H,~a\,E_H+a\,\Delta_H\Big]. \tag{2.20}
\end{align*}

\section{Basic Properties of the Continuous Special Affine Wavelet Transform}

\parindent=0mm \vspace{.0in}
In this Section, we shall study some fundamental properties of the proposed special affine wavelet transform (2.4) such as Rayleigh’s theorem, inversion formula  and characterization of the range are discussed using the machinery of special affine Fourier transforms and operator theory. In this direction, we have the following theorem which assembles some of the basic properties of the proposed transform.

\parindent=0mm \vspace{.1in}
{\bf Theorem 3.1.} {\it For any $f,g\in L^2(\mathbb R)$  and $\alpha,\beta\in\mathbb R$,  the  continuous special affine wavelet transform  defined by  (2.4)  satisfies the following properties:}

\parindent=0mm \vspace{.1in}
(i)~~${\text {\it Linearity:}}~~\qquad{\mathcal{W}}^{M}_\psi \Big[\alpha f(t)+\beta g(t)\Big] (a,b)= \alpha~{\mathcal{W}}^{M}_\psi \big[f\big](b, a)+\beta~{\mathcal{W}}^{M}_\psi \big[g\big](a,b)$.

\parindent=0mm \vspace{.1in}
(ii)~${\text {\it Translation:}} {\mathcal{W}}^{M}_\psi \Big[f(t-\alpha)\Big](a,b)=\exp\left\{\dfrac{i}{2B}\Big(A\alpha^2-\alpha(p-a)\Big)\right\} {\mathcal{W}}^{M}_\psi\left[\exp\left\{\dfrac{iAt\alpha}{B}\right\} f(t)\right] (b-\alpha,a).$

\parindent=0mm \vspace{.1in}
(iii)~~{\text {\it Parity:}}~$ {\mathcal{W}}^{M}_\psi\Big[f(-t)\Big](a,b)= (-i)\, {\mathcal{W}}^{M^\prime}_\psi\big[f\big](-b,-a),\quad \text{where}\quad M^\prime =(A, B, C, D, -p, -q) $.

\parindent=0mm \vspace{.1in}
(iv)~~{\text {\it Dilation:}}~\quad  ${\mathcal{W}}^{M}_\psi\Big[D_\alpha f\Big](a,b)={\mathcal{W}}^{M^\prime}_\psi\big[f\big]\left(\dfrac{a}{\alpha},\dfrac{b}{\alpha}\right),  \quad\text{where}\quad M^\prime=\big(aA,aB,aC,aD,p,q\big)$.

\parindent=0mm \vspace{.1in}
(v)~${\text {\it Conjugation:}}~~\quad {\mathcal{W}}^M_\psi\Big[\overline {f}\Big] (b,a)=\overline{{\mathcal{W}}^{M^{-1}}_{\overline {\psi}}\big[f\big](b,a)}.$

\parindent=0mm \vspace{.1in}

{\it Proof.} For the sake of brevity, we omit the proof of (i).

\parindent=0mm \vspace{.1in}
(ii). For any  $\alpha\in \mathbb R$, we have
\begin{align*}
&{\mathcal{W}}^{M}_\psi \Big[f(t-\alpha)\Big](a,b)\\
&=\dfrac{1}{\sqrt{2\pi iB}}\int_{-\infty}^{\infty} f(t-\alpha)\dfrac{1}{\sqrt{a}} \overline{\psi\left(\dfrac{t-b}{a}\right)} \exp\left\{\dfrac{i}{2B}\Big(At^2+2t(p-a)) \right\}\\
&\qquad\qquad\qquad\qquad\qquad\qquad\qquad\times \exp\left\{\dfrac{i}{2B}\Big(-2a(Dp-Bq)+D(a^2+p^2)\Big) \right\}dt\\
&=\dfrac{1}{\sqrt{2\pi iBa}} \int_{-\infty}^{\infty} f(t-\alpha)\overline{\psi\left(\dfrac{t-\alpha-(b-\alpha)}{a}\right)}\exp\left\{\dfrac{i}{2B}\Big(A(t-\alpha)^2-A\alpha^2+2A\alpha t\Big) \right\}\\
&\qquad\qquad \times \exp\left\{\dfrac{i}{2B}\Big(2(t-\alpha)(p-a)+\alpha(p-a)-2a(Dp-Bq)+D(a^2+p^2)\Big) \right\}dt\\
&=\exp\left\{\dfrac{i}{2B}\Big(-A\alpha^2+\alpha(p-a)\Big) \right\} \dfrac{1}{\sqrt{2\pi iBa}} \int_{-\infty}^{\infty} f(t-\alpha)\exp\left\{\dfrac{i\alpha At}{B}\right\}\, \overline{\psi_{a,b-\alpha}(t-\alpha)} \mathcal{K}^M(t-\alpha,a)dt\\
&=\exp\left\{\dfrac{i}{2B}\Big(A\alpha^2-\alpha(p-a)\Big)\right\} {\mathcal{W}}^{M}_\psi\left[\exp\left\{\dfrac{iAt\alpha}{B}\right\}\big[f\big]\right](b-\alpha,a).
\end{align*}

(iii). To check the behaviour of the special affine wavelet transform \eqref{2.3} under the reflection of the input signal $f$, we proceed as

\begin{align*}
&{\mathcal{W}}^{M}_\psi\Big[f(-t)\Big](a,b)\\
&= \int_{-\infty}^{\infty} f(-t)\dfrac{1}{\sqrt{a}} \overline{\psi\left(\dfrac{(t-b)}{a}\right)}\, \dfrac{1}{\sqrt{2\pi iB}} \exp\left\{\dfrac{i}{2B}\Big(At^2+2t(p-a)\Big) \right\}\\
&\qquad\qquad\qquad\qquad\qquad\qquad\times \exp\left\{\dfrac{i}{2B}\Big(-2a(Dp-Bq)+D\big(a^2+p^2\big)\Big) \right\}dt\\
&= \int_{-\infty}^{\infty} f(-t)\dfrac{1}{\sqrt{a}} \overline{\psi\left(\dfrac{-t-(-b)}{-a}\right)}\, \dfrac{1}{\sqrt{2\pi iB}} \exp\left\{\dfrac{i}{2B}\Big(A(-t)^2+2(-t)\big(-p-(-a)\big)\Big) \right\}\\
&\qquad\qquad\qquad\times \exp\left\{\dfrac{i}{2B}\Big(-2(-a)\Big(D(-p)-B(-q)\Big)+D\big((-a)^2+(-p)^2\big)\Big) \right\}dt
\\&= \int_{-\infty}^{\infty} f(t^\prime)\dfrac{i}{\sqrt{-a}} \overline{\psi\left(\dfrac{t^\prime-(-b)}{-a}\right)}\, \dfrac{1}{\sqrt{2\pi iB}} \exp\left\{\dfrac{i}{2B}\Big(A(t^\prime)^2+2(t')\big(-p-(-a)\big)\Big) \right\}\\
&\qquad\qquad\qquad\times \exp\left\{\dfrac{i}{2B}\Big(-2(-a)\big(D(-p)-B(-q)\big)+D\big((-a)^2+(-p)^2\big)\Big) \right\}(-dt^\prime)\\
&= (-i)\, {\mathcal{W}}^{M^\prime}_\psi\Big[f\Big](-b,-a).
\end{align*}

(iv). For any $\alpha\in\mathbb R$, we have
\begin{align*}
&{\mathcal{W}}^{M}_\psi\big[D_\alpha f\big](a,b)\\
&= \int_{-\infty}^{\infty}\dfrac{1}{\sqrt{\alpha}} f\left(\dfrac{t}{\alpha}\right)\, \dfrac{1}{\sqrt{a}} \overline{\psi\left(\dfrac{t-b}{a}\right)}\, \dfrac{1}{\sqrt{2\pi iB}} \exp\left\{\dfrac{i}{2B}\Big(At^2+2t(p-w)\Big) \right\}\\
&\qquad\qquad\qquad\qquad\qquad\qquad\times \exp\left\{\dfrac{i}{2B}\Big(-2w(Dp-Bq)+D(w^2+p^2)\Big) \right\}dt\\
&= \int_{-\infty}^{\infty}\left(\dfrac{\sqrt{\alpha}}{\sqrt{a}}\right)^2 f\left(t^\prime\right) \overline{\psi\left(\dfrac{\alpha t^\prime-b}{a}\right)} \dfrac{1}{\sqrt{2\pi iB}} \exp\left\{\dfrac{i}{2B}\Big(A\big(\alpha t^\prime\big)^2+2\alpha t(p-w)\Big) \right\}\\
&\qquad\qquad\qquad\qquad\qquad\qquad\times \exp\left\{\dfrac{i}{2B}\Big(-2w(Dp-Bq)+D(w^2+p^2)\Big) \right\}dt\\
&= \left(\sqrt{\dfrac{\alpha}{a}}\right)\,\int_{-\infty}^{\infty} f\left(t^\prime\right) \overline{\psi\left(\dfrac{\alpha \Big( t^\prime-\frac{b}{a}\Big)}{a}\right)} \dfrac{1}{\sqrt{2\pi iB^\prime}} \exp\left\{\dfrac{i}{2B^\prime}\Big(A^\prime\big( t^\prime\big)^2+2 t^\prime(p-w)\Big) \right\}\\
&\qquad\qquad\qquad\qquad\qquad\qquad\times \exp\left\{\dfrac{i}{2B^\prime}\Big(-2w(D^\prime p-B^\prime q)+D^\prime(w^2+p^2)\Big) \right\}dt^\prime\\
&={\mathcal{W}}^{M^\prime}_\psi\Big[f\Big]\left(\dfrac{a}{\alpha},\dfrac{b}{\alpha}\right).
\end{align*}

\parindent=0mm \vspace{.0in}
(v). It is easy to compute  and is therefore omitted.

\parindent=0mm \vspace{.1in}
This completes the proof of Theorem 3.1. \quad\fbox

\parindent=8mm \vspace{.1in}

 In the next theorem, our aim is to establish an orthogonality formula between two signals and their respective special affine wavelet transforms. As a consequence of this formula, we can deduce the resolution of identity for the proposed transform (2.4).

 \parindent=0mm \vspace{.1in}
{\bf Theorem 3.2 {(\bf Rayleigh's Theorem).}}  {\it Let ${\mathcal{W}}_\psi^M \big[f\big](a,b)$  and $ {\mathcal{W}}_\phi^M \big[g\big](a,b)$ be the special affine wavelet transforms of  $f$ and $g$ belonging to $L^2(\mathbb R)$, respectively. Then, we have}
\begin{align*}\label{3.1}
\Big\langle {\mathcal{W}}_\psi^M \big[f\big](a,b),{\mathcal{W}}_\phi^M \big[g\big](a,b) \Big\rangle = \Big\langle f,g \Big\rangle\overline{\Big\langle \psi_{a,b}^M,\phi_{a,b}^M \Big\rangle}.\tag{3.1}
\end{align*}

 \parindent=0mm \vspace{.1in}
{\it Proof.} Since
\begin{align*}
\int_{-\infty}^{\infty} {\mathcal K}^{M}(t,a)\, \overline{{\mathcal K}^{M}(x,a)}\, da = \delta(t-x)
\end{align*}
therefore, we have
\begin{align*}
&\Big\langle {\mathcal{W}}_\psi^M \big[f\big](a,b),{\mathcal{W}}_\phi^M \big[g\big](a,b) \Big\rangle \\ &\qquad\qquad=\int_{-\infty}^{\infty}\int_{-\infty}^{\infty}{\mathcal{W}}_\psi^M \big[f\big](a,b) \, \overline{ {\mathcal{W}}_\psi^M \big[g\big](a,b)}\,\, dadb\\
&\qquad\qquad=\int_{-\infty}^{\infty}\int_{-\infty}^{\infty}\left\{\int_{-\infty}^{\infty} f(t)\, \overline{\psi_{a,b}^M}(t)\, \mathcal{K}^M(t,a)\,dt \right\}\overline{\left\{\int_{-\infty}^{\infty} g(t^\prime)\, \overline{\phi_{a,b}^M}(t^\prime)\, \mathcal{K}^M(t^\prime,a)\,dt^\prime\right\}}\\
&\qquad\qquad=\int_{-\infty}^{\infty}\int_{-\infty}^{\infty}\left\{\int_{-\infty}^{\infty}\int_{-\infty}^{\infty}f(t)\overline{g(t^\prime)}\,\overline{\psi_{a,b}^M(t^\prime)}\phi_{a,b}^M(t) \overline {\mathcal{K}^M(t,a)} \mathcal{K}^M(t^\prime,a)\,dtdt^\prime \right\}\,\,dadb\\
&\qquad\qquad= \int_{-\infty}^{\infty}\left\{\int_{-\infty}^{\infty}\int_{-\infty}^{\infty}f(t)\overline{g(t)}\,\overline{\psi_{a,b}^M(t)}\phi_{a,b}^M(t)\,dt \right\}\,dadb\\
&\qquad\qquad= \int_{-\infty}^{\infty}f(t)\overline{g(t)}\,\bigg\{\int_{-\infty}^{\infty}\int_{-\infty}^{\infty}\overline{\psi_{a,b}^M(t)\, \phi_{a,b}^M(t)} \,\, dadb\bigg\}\,dt\\
&\qquad\qquad= \Big\langle f,g \Big\rangle\overline{\Big\langle \psi_{a,b}^M,\phi_{a,b}^M\Big\rangle}.
\end{align*}
This completes the proof of Theorem 3.2.

\parindent=0mm \vspace{.1in}
{\it Remarks:} Theorem 3.2 allows us to make the following remarks:

\parindent=0mm \vspace{.1in}
(i) For $f=g$, the identity (3.1) boils down to Moyal's principle
\begin{align*}
\int_{-\infty}^{\infty}\int_{-\infty}^{\infty}\Big|{\mathcal{W}}_\psi^M \big[f\big](b,a)\Big|^2da\,db =\big\|f\big\|^2_2\, \big\|\psi\big\|^2_2.\tag{3.2}
\end{align*}
(ii) For a normalized window function  $\psi\in L^2(\mathbb R)$, we have
\begin{align*}
\int_{-\infty}^{\infty}\int_{-\infty}^{\infty}\Big|{\mathcal{W}}_\psi^M \big[f\big](b,a)\Big|^2da\,db =\big\|f \big\|^2_2. \tag{3.3}
\end{align*}
(iii) For the normalized functions $f,\psi\in L^2(\mathbb R)$, equation (3.2) gives the Radar uncertainty principle for the special affine wavelet transform
\begin{align*}
\int_{-\infty}^{\infty}\int_{-\infty}^{\infty}\Big|{\mathcal{W}}_\psi^M \big[f\big](b,a)\Big|^2 da\,db = 1.  \tag{3.4}
\end{align*}

\parindent=0mm \vspace{.0in}
After the complete formulation of the special affine wavelet transform (2.4), our next goal is to find a way out for the reconstruction of the input signal $f\in L^2(\mathbb R)$  via the coefficients of $\left\langle f, \psi_{a,b}^M\right\rangle$.

\parindent=0mm \vspace{.1in}
{\bf Theorem 3.3 ( Inversion Formula).} {\it If $~{\mathcal{W}}_\psi^M \big[f\big](a,b)$ is the special affine wavelet transform of an arbitrary function $f\in L^2(\mathbb R)$, then $f$ can be reconstructed as}
\begin{align*}\label{3.5}
f(t)= \dfrac{1}{\overline{\big\langle \psi_{a,b}^M,\phi_{a,b}^M \big\rangle}} \int_{-\infty}^{\infty}\int_{-\infty}^{\infty}{\mathcal{W}}_\psi^M \big[f\big](a,b) \, \overline{\mathcal{K}^M(t,a)} \, \phi_{a,b}^M(t)\,da\,db, \quad a.e. \tag{3.5}
\end{align*}

\parindent=0mm \vspace{.1in}
{\it Proof.} Implication of the Moyal's principle (3.2) yields that
\begin{align*}
g(t)= \dfrac{1}{\overline{\big\langle \psi_{a,b}^M,\phi_{a,b}^M \big\rangle}} \int_{-\infty}^{\infty}\int_{-\infty}^{\infty}{\mathcal{W}}_\psi^M \big[f\big](a,b) \, \overline{\mathcal{K}^M(t,a)} \, \phi_{a,b}^M(t)\,da\,db
\end{align*}
is well defined in $L^2(\mathbb R)$. Moreover, for any $g,h\in L^2(\mathbb R)$, the Rayleigh's Theorem 3.2 gives
\begin{align*}
\big\langle g,\,h \big\rangle &= \dfrac{1}{\overline{\big\langle \psi_{a,b}^M,\phi_{a,b}^M \big\rangle}} \int_{-\infty}^{\infty}\,\left\{\int_{-\infty}^{\infty}\int_{-\infty}^{\infty}{\mathcal{W}}_\psi^M \big[f\big](a,b) \, \overline{\mathcal{K}^M(t,a)} \phi_{a,b}(t)\,da\,db \, \overline{h(t)}\right\}\, dt \\
&=\dfrac{1}{\overline{\big\langle \psi_{a,b}^M,\phi_{a,b}^M \big\rangle}} \int_{-\infty}^{\infty} \int_{-\infty}^{\infty}\int_{-\infty}^{\infty}{\mathcal{W}}_\psi^M \big[f\big](a,b) \, \overline{h(t) \mathcal{K}^M(t,a)} \phi_{a,b}(t)\,\,da\,db\,dt \\
&=\dfrac{1}{\overline{\big\langle \psi_{a,b}^M,\phi_{a,b}^M \big\rangle}} \int_{-\infty}^{\infty}\int_{-\infty}^{\infty}{\mathcal{W}}_\psi^M \big[f\big](a,b) \, \left\{ \int_{-\infty}^{\infty}\overline{h(t)}\,\overline{ \mathcal{K}^M(t,a)} \phi_{a,b}(t)\, dt \right\}\,\,da\,db \\
&=\dfrac{1}{\overline{\big\langle \psi_{a,b}^M,\phi_{a,b}^M \big\rangle}} \int_{-\infty}^{\infty}\int_{-\infty}^{\infty}{\mathcal{W}}_\psi^M \big[f\big](a,b)\, \overline{{\mathcal{W}}_\psi^M \big[f\big](a,b)} \,\,da\,db \\
&=\dfrac{1}{\overline{\big\langle \psi_{a,b}^M,\phi_{a,b}^M \big\rangle}}\Big\langle {\mathcal{W}}_\psi^M \big[f\big](a,b),{\mathcal{W}}_\psi^M \big[h\big](a,b) \Big\rangle\\
&=\Big\langle f,\,h\Big\rangle.
\end{align*}
Since $h$ is  an arbitrary element of $L^2(\mathbb R)$, therefore it follows that
\begin{align*}
f(t)= \dfrac{1}{\overline{\big\langle \psi_{a,b}^M,\phi_{a,b}^M \big\rangle}} \int_{-\infty}^{\infty}\int_{-\infty}^{\infty}{\mathcal{W}}_\psi^M \big[f\big](a,b) \, \overline{\mathcal{K}^M(t,a)} \, \phi_{a,b}^M(t)\,da\,db, \quad a.e.
\end{align*}
This completes the proof of Theorem 3.3. \qquad \fbox

\parindent=0mm \vspace{.1in}
{\it Remark.} For the normalized functions $\psi, \phi \in L^2(\mathbb R)$ with $\psi=\phi$,  equation (3.5) reduces to a relatively compact  formula
\begin{align*}
f(t)=  \int_{-\infty}^{\infty}\int_{-\infty}^{\infty}{\mathcal{W}}_\psi^M \big[f\big](a,b) \, \overline{\mathcal{K}^M(t,a)} \, \psi_{a,b}^M(t)\,da\,db, \quad a.e.\tag{3.6}
\end{align*}

\parindent=0mm \vspace{.0in}

In the upcoming theorem, we  shall obtain a complete characterization of the range of the transform ${\mathcal{W}}^{M}_{\psi}\big[f\big](a,b)$ defined in (2.4). The result follows as a consequence of the  inversion formula (3.5) and the well known Fubini theorem.

\parindent=0mm \vspace{.1in}
{\bf Theorem 3.4 (Characterization of Range).} {\it If $ h \in L^2(\mathbb {R}^2)$ and $\psi$ is a normalized wavelet. Then, $h$ is the special affine wavelet transform of a certain square integrable function if and only if it satisfies the reproduction property:}
\begin{align*}\label{3.6}
h(c,d) = \int_{-\infty}^{\infty}\int_{-\infty}^{\infty} h(a,b)\, \Big\langle \psi^M_{a,b},\,\psi^M_{c,d}\Big\rangle\,da\,db.\tag{3.7}
\end{align*}

\parindent=0mm \vspace{.0in}
{\it Proof.} Let $h$ belongs to the range of the transform ${\mathcal{W}}^{M}_{\psi}$. Then, there exists a function $f\in L^2(\mathbb R)$, such that ${\mathcal{W}}^{M}_{\psi} \big[f\big] = h$. In order to show that $h$ satisfies (3.7), we proceed as
\begin{align*}
h(c,d)&= {\mathcal{W}}^{M}_{\psi}\big[f\big](c,d)\\
&=\int_{-\infty}^{\infty} f(t)\, \overline{\psi^M_{c,d}(t)}\,dt\\
&=\int_{-\infty}^{\infty} \dfrac{1}{\big\|\psi \big\|^2_2} \left\{\int_{-\infty}^{\infty}\int_{-\infty}^{\infty}{\mathcal{W}}^{M}_{\psi}\big[f\big] (a,b) \, \psi^M_{a,b}(t) \,dadb \right\} \overline{\psi^M_{c,d}(t)}\,dt\\
&=\int_{-\infty}^{\infty} \int_{-\infty}^{\infty} {\mathcal{W}}^{M}_\psi \big[f\big] (a,b) \left\{ \int_{-\infty}^{\infty} \psi^M_{a,b}(t)\,\overline{\psi^M_{c,d}(t)} dt \right\} \,da\,db  \\
&=\int_{-\infty}^{\infty} \int_{-\infty}^{\infty} h(a,b)\, \Big\langle \psi^M_{a,b},{\psi^M_{c,d}}\Big\rangle\, da\,db.
\end{align*}

Conversely, suppose that an arbitrary function $h\in L^2(\mathbb {R}^2)$ satisfies (3.7). Then, we show that there exist $f\in L^2(\mathbb R)$, such that ${\mathcal{W}}^{M}_{\psi}\big[f\big] = h$. Assume that
\begin{align*}
f(t)=\int_{-\infty}^{\infty} \int_{-\infty}^{\infty} h(a,b)\, \psi^M_{a,b}(t)\, da\,db.
\end{align*}
Firstly, we claim that $f\in L^2(\mathbb R).$ To prove this, we proceed as
\begin{align*}
\big\|f\big\|^2_2&=\int_{-\infty}^{\infty} f(t)\,\overline{f(t)}\,dt\\
&=\int_{-\infty}^{\infty} \left\{ \int_{-\infty}^{\infty} \int_{-\infty}^{\infty} h(a,b)\, \psi^M_{a,b}(t)\, dadb \right\}
 \left\{ \int_{-\infty}^{\infty} \int_{-\infty}^{\infty} \overline{ h(a,b)\, \psi^M_{a,b}(t)}\, dadb \right\} \, dt\\
&=\int_{-\infty}^{\infty} \int_{-\infty}^{\infty} \int_{-\infty}^{\infty} \int_{-\infty}^{\infty} h(a,b) \overline{ h(a,b)}\left\{ \int_{-\infty}^{\infty} \psi^M_{a,b}(t)\overline{\psi^M_{a,b}(t)}\,dt\right\}\,dadb\,dadb\\
&=\int_{-\infty}^{\infty} \int_{-\infty}^{\infty} \overline{ h(a,b)}\,\left\{\int_{-\infty}^{\infty} \int_{-\infty}^{\infty} h(a,b)\big\langle \psi^M_{a,b},\psi^M_{c,d}\big\rangle \, dadb \right\}\, dadb\\
&=\int_{-\infty}^{\infty} \int_{-\infty}^{\infty} \overline{ h(a,b)}\,h(a,b)\,dadb\\
&=\int_{-\infty}^{\infty} \int_{-\infty}^{\infty} \big\|h(a,b)\big\|^2_2 da\,db\\
&= \big\|h\big\|^2_2,
\end{align*}
which establishes our claim . Moreover, as a consequence of the well-known Fubini-theorem, we have
\begin{align*}
{\mathcal{W}}^{M}_{\psi}\big[f\big](c,d)&= \int_{-\infty}^{\infty} f(t)\, \overline{\psi^M_{c,d}(t)}\, dt\\
&= \int_{-\infty}^{\infty} \left\{\int_{-\infty}^{\infty} \int_{-\infty}^{\infty} h(a,b)\, \psi^M_{a,b}(t)\, dadb\right\} \, \overline{\psi^M_{c,d}(t)}\, dt\\
&= \int_{-\infty}^{\infty} \int_{-\infty}^{\infty} h(a,b) \, \int_{-\infty}^{\infty} \psi^M_{a,b}(t) \overline{\psi^M_{c,d}(t)} dt\, da\, db \\
&= \int_{-\infty}^{\infty} \int_{-\infty}^{\infty} h(a,b) \, \Big\langle \psi^M_{a,b}, \psi^M_{c,d}\Big\rangle\, da\,db \\
&= h(c,d).
\end{align*}
This evidently completes the proof of Theorem 3.4. \qquad\fbox

\parindent=0mm \vspace{.1in}

{\bf Corollary 3.5 (Reproducing Kernel Hilbert space).}  {\it For a normalized wavelet $\psi\in L^2(\mathbb R)$, the range of the continuous special affine wavelet transform (2.4) is a reproducing kernel Hilbert space in $L^2(\mathbb R^2)$ with kernel given by}
\begin{align*}
\mathcal{K}^M_\psi(a,b,c,d) = \Big\langle \psi^M_{a,b},\, \psi^M_{c,d}\Big\rangle.
\end{align*}

\section{Discrete Special Affine Wavelet Transform}

\parindent=0mm \vspace{.0in}
In this Section, we introduce a discrete analogue of the proposed special affine wavelet transform (2.4) and obtained the associated reconstruction formula. Our strategy is based on an appropriate discretization of the parameters $a$ and $b$ involved in the newly introduced family $\psi^M_{a, b}(t)$ given by (2.2). To facilitate this, we choose $a = a_0^j $  and $ b =kb_0a_0^j,$, where $j,k \in \mathbb Z$ and $a_0,\,b_0$ are fixed positive constants. Consequently, the desired discretized family takes the following form:
\begin{align*}\label{4.1}
\psi^M_{j,k}(t)= a_0^{-j/2} \psi\left( a_0^{-j}t-kb_0\right)\, {\mathcal K}^M\left(t,kb_0a_0^{-j}\right),\tag{4.1}
\end{align*}
where ${\mathcal K}^M(t,kb_0a_0^{-j})$ is the discrete version of (2.3). Formally, we have the following definition of the discrete special affine wavelet transform.

\parindent=0mm \vspace{.1in}
{\bf Definition 4.1.} For any arbitrary function $f\in L^2(\mathbb R)$ and uni-modular matrix $M=(A,B,C,D,p,q)$, the discrete special affine wavelet transform with respect  the window function $\psi \in L^2(\mathbb{R})$ is given by
\begin{align*}\label{4.2}
{\mathcal{W}}_\psi^M (j,k)= \int_{-\infty}^{\infty} f(t)\,\overline{\psi^{M}_{j,k}(t)}\,dt,\tag{4.2}
\end{align*}
where $\psi^M_{j,k}(t)$ is given by (4.1).

\parindent=8mm \vspace{.1in}
For the illustration of the discrete special affine wavelet transform (4.2), we have the following example.

\parindent=0mm \vspace{.1in}
{\bf Example 4.2.} We shall recall the generalized Gaussian function $f(t)=e^{-iAt^2/2B},~AB>0$  and from Example 2.2(ii), and compute its discrete special affine wavelet transform with respect to the window function
\begin{align*}
\psi(t)=\left\{\begin{array}{cc}
1, & 0\le t\le 1/2 \\
-1, & 1/2 \le t\le 1 .
\end{array}\right.\tag{4.3}
\end{align*}
The dyadic dilation and integer translations of $\psi(t)$ are given by
\begin{align*}
\psi\left( 2^{-j}t-k \right) =\left\{\begin{array}{cc}
1, & 2^jk  \le t\ < 2^j \left(k+ \frac{1}{2}\right) \\
-1, & 2^j\left(k+\frac{1}{2}\right) \le t\le 2^j(k+1).
\end{array}\right.\tag{4.4}
\end{align*}
Therefore, the discrete special affine wavelet transform of $f(t)$ can be computed as
\begin{align*}
&{\mathcal{W}}^M_\psi(j,k)\quad\\
&= \int_{-\infty}^{\infty} \exp\left\{{\frac{-iAt^2}{2B}}\right\} 2^{-j/2} \,\overline{\psi\left(2^{-j/2}t-k\right)}\\
&\quad \times \exp\left\{\dfrac{i}{2B}\Big(At^2+2t(p-k2^j)-2^{j+1}k(Dp-Bq)+D\big(k^22^{2j}+p^2\big) \Big)\right\}\,dt\\
&=2^{-j/2}\,\exp\left\{\dfrac{i}{2B}\Big(-2^{j+1}k(Dp-Bq)+D\big(k^2 2^{2j}+p^2\big) \Big)\right\}\\
&\quad \times \Bigg[ \int_{2^jk}^{2^j\left(k+\frac{1}{2}\right)}\exp\left\{\dfrac{i}{B}\big(p-k2^j\big)t\right\}\, dt +  \int_{2^j\left(k+\frac{1}{2}\right)}^{2^j(k+1)}\,(-1)\,\exp\left\{\dfrac{i}{B}\big(p-k2^j\big)t\right\}\, dt \, \Bigg]\\
&=2^{-j/2}\,\exp\left\{\dfrac{i}{2B}\Big(-2^{j+1}k(Dp-Bq)+D\big(k^2 2^{2j}+p^2\big) \Big)\right\}\\
&\quad \times \dfrac{B}{i\big(p-k2^j\big)} \Bigg[ \exp\left\{\dfrac{i}{B}\big(p-k2^j\big)2^j\left(k+\frac{1}{2}\right) \right\} - \exp\left\{\dfrac{i}{B}\big(p-k2^j\big)2^jk\right\}\\
&\qquad -\exp\left\{\dfrac{i}{B}\big(p-k2^j\big)2^j(k+1) \right\} + \exp\left\{\dfrac{i}{B}\big(p-k2^j\big)2^j\left(k+\frac{1}{2}\right) \right\} \Bigg]\\
&= \dfrac{B\,2^{-j/2}}{i(p-k2^j)} \,\exp\left\{\dfrac{i}{2B}\Big(-2^{j+1}k(Dp-Bq)+D\big(k^2 2^{2j}+p^2\big) \Big)\right\}\\
&\quad \times \exp\left\{\dfrac{i}{B}\big(p-k2^j\big)2^jk \right\}\Bigg[2\exp\left\{\dfrac{i}{2B}\big(p-k2^j\big) -1 \right\} -\exp\left\{\dfrac{i}{B}\big(p-k2^j\big)2^j \right\}  \Bigg].
\end{align*}

\parindent=0mm \vspace{.0in}
We are now focussed to establish a reconstruction formula for the discrete special affine wavelet transform defined by \eqref{4.2}.

\parindent=0mm \vspace{.1in}
{\bf Theorem 4.3.} {\it Let ${\mathcal{W}}^M_\psi(j,k)$ discrete special affine wavelet transform of $f\in L^2(\mathbb R)$, defined by (4.2). Assume that the discrete system $\big\{\psi^M_{j,k}; j,k\in \mathbb Z\big\}$ as defined by (4.1) satisfies the stability condition}
\begin{align*}\label{4.3}
E\big\|f\big\|_2^2\le\sum_{j\in\mathbb{Z}}\sum_{k\in\mathbb{Z}}\Big|\Big\langle f, \psi^M_{j,k}\Big\rangle_2\Big|^2\le F\big\|f\big\|_2^2,\tag{4.5}
\end{align*}
and let $S$ be the given linear operator
\begin{align*}\label{4.4}
 S(f)= \sum_{j\in\mathbb{Z}}\sum_{k\in\mathbb{Z}}\Big\langle f,\,\psi^M_{j,k} \Big\rangle \psi^M_{j,k}. \tag{4.6}
\end{align*}
Then, $f(t)$ can be reconstructed as
\begin{align*}
 f(t)= \sum_{j\in\mathbb{Z}}\sum_{k\in\mathbb{Z}} \Big\langle f,\,\psi^M_{j,k} \Big\rangle S^{-1} \left(\psi^M_{j,k}\right). \tag{4.7}
\end{align*}

{\it Proof.} By virtue of the stability condition \eqref{4.3}, it can be easily verified that the operator $S$ defined by \eqref{4.4} is one-one and bounded. Then, for $S(f)=h$, we have
\begin{align*}
\big\langle h,f\big\rangle&=\Big\langle S(f),f\Big\rangle\\
&=\Big\langle \sum_{j\in\mathbb{Z}}\sum_{k\in\mathbb{Z}}\Big\langle f, \psi^M_{j,k}\Big\rangle \psi^M_{j,k},\,f\Big\rangle\\
&=\sum_{j\in\mathbb{Z}}\sum_{k\in\mathbb{Z}}\Big\langle f, \psi^M_{j,k}\Big\rangle \overline{\Big\langle f, \psi^M_{j,k}\Big\rangle}\\
&=\sum_{j\in\mathbb{Z}}\sum_{k\in\mathbb{Z}}\Big|\Big\langle f, \psi^M_{j,k}\Big\rangle\Big|^2
\\&\le F\big\|f\big\|_2^2.
\end{align*}

Moreover, an implication of Cauchy-Schwartz inequality yields
\begin{align*}
E\,\left\|S^{-1}(h)\right\|^2_2&=E\,\left\|S^{-1}\big(S(f)\big)\right|^2_2\\
&=E\,\big\|f\big\|^2_2\\
&\le \sum_{j\in\mathbb{Z}}\sum_{k\in\mathbb{Z}}\Big|\Big\langle f, \psi^M_{j,k}\Big\rangle\Big|^2\\
&=\big\langle h,f\big\rangle\\
&\le\big\|h\big\|_2\,\big\|f\big\|_2\\
&=\big\|h\big\|_2\,\big\|S^{-1}(h)\big\|_2.
\end{align*}
which implies that $E\,\big\|S^{-1}(h)\big\|_2\le\big\|h\big\|_2$. Therefore, for every $f\in L^2(\mathbb R)$, we have
\begin{align*}
f=S^{-1}S(f)=S^{-1}\bigg\{\sum_{j\in\mathbb{Z}}\sum_{k\in\mathbb{Z}}\Big\langle f, \psi^M_{j,k}\Big\rangle \psi^M_{j,k}\bigg\}=\sum_{j\in\mathbb{Z}}\sum_{k\in\mathbb{Z}}\Big\langle f, \psi^M_{j,k}\Big\rangle\,S^{-1}\left(\psi^M_{j,k}\right).
\end{align*}
This completes the proof of the Theorem 4.3. \qquad\fbox

\section{Relationship Between The Special Affine Wigner Distribution and The Proposed Wavelet Transform}

\parindent=0mm \vspace{.0in}
The Wigner-Ville distribution (WVD) is defined as the Fourier transform of the instantaneous autocorrelation function $f\big(t+\frac{\tau}{2}\big)\overline{f\big(t-\frac{\tau}{2}\big)}$ and is regarded as the main distribution of all the time-frequency distributions. It is the premier tools in time-frequency analysis with applications in diverse areas including optical and radar systems,  and non- stationary signal processing in general. The vast applicability of the WVD in time-frequency analysis has prompted serious attention in recent years (See \cite{DS}). In this direction,  Urynbassarova et al.\cite{ULT} introduced a generalised Wigner-Ville distribution namely, the special affine Wigner distribution by substituting the Fourier transform kernel with the SAFT kernel and verified that the new Wigner distribution is more effective in applications that both the linear canonical Wigner distribution as well as the classical WVD. In this Section, we shall develop a relationship between the special affine wavelet transform as defined by \eqref{2.3}  and the special affine Wigner distribution.

\parindent=0mm \vspace{.1in}

{\bf Definition 5.1.} For any $f\in L^2(\mathbb{R})$, the special affine Wigner distribution with uni-modular matrix $M =(A,B,C,D,p,q)$  is defined by \cite{}
\begin{align*}\label{5.1}
\Big[{WD}^M f\Big](t,a)=\int_{-\infty}^{\infty} f\left(t+\dfrac{\tau}{2}\right)\overline{f\left(t-\frac{\tau}{2}\right)}\, \mathcal{K}^M(\tau,a)\,d\tau,\tag{5.1}
\end{align*}
where $\mathcal{K}^M(\tau,a)$ is given by (2.3).

\parindent=8mm \vspace{.1in}
In the following theorem, we  derive a relationship between special affine Wingner distribution and the proposed transform \eqref{2.3}.

\parindent=0mm \vspace{.1in}
{\bf Theorem 5.2.} {\it If $\psi\in L^2(\mathbb{R})$, ${\mathcal{W}}^M_\psi\big[f\big](a,b)$ and  ${WD}^M \big[f\big](t,a)$ are the continuous especial affine wavelet transform and special affine Wigner distribution of $f\in L^2(\mathbb R)$ defined by (4.2)  and (5.1), respectively. Then,  we have}
\begin{align*}
\Big[{WD}^M f\Big](t,a) &= \dfrac{1}{2\sqrt {2\pi iB}}\, \dfrac{1}{\big\|\psi\big\|^2_2}\, \exp\left\{\dfrac{i}{2B}\Big(-2a(Du_{0}-Dq)+D\big(a^2+u^2_{0}\big)\Big) \right\}\\
& \times \int_{-\infty}^{\infty}\int_{-\infty}^{\infty}{\mathcal{W}}^M_\psi \left[ \exp\left\{\dfrac{i}{B}\Big(2Az^2+2z(u_{0}-a-2At)\Big) f(t)\right\}\right](b,a)\, dz\\
&\times \overline{{\mathcal{W}}^M_\psi \left[\exp\left\{\dfrac{-2At^2}{B}\right\} f(t)\right]} (-b-2t,-a) \,da\,db.\tag{5.2}.
\end{align*}

{\it Proof.} By making a change of variables as  $t+\frac{\tau}{2}=z$ in \eqref{5.1}, we obtain
\begin{align*}
&\Big[{WD}^M f\Big](t,a)\\
&=\dfrac{1}{\sqrt {2\pi iB}}\exp\left\{\dfrac{i}{2B}\Big(-2a(Du_{0}-Dq)+D\big(a^2+u_{0}^2\big)\Big) \right\}\\
&\quad\times \int_{-\infty}^{\infty} f(z)\, \overline{f(2t-z)}\cdot\exp\left\{\dfrac{i}{2B}\Big(A(2z-2t)^2+2(2z-2t)(p-a)\Big) \right\}\, \dfrac{dz}{2}\\
&=\dfrac{1}{\sqrt {2\pi iB}}\exp\left\{\dfrac{i}{2B}\Big(-2a(Du_{0}-Dq)+D\big(a^2+u_{0}^2\big)\Big) \right\}\\
&\quad \times \int_{-\infty}^{\infty} f(z)\, \overline{f(2t-z)}\,\exp\left\{\dfrac{i}{B}\Big(2A(z^2+t^2-2zt)^2+2(z-t)(p-a)\Big) \right\}\,dz. \tag{5.3}
\end{align*}
Exploiting  the inversion formula (3.5) and translation property of the continuous special affine wavelet transform (2.4), we get

\begin{align*}
f(2t-z)&=\dfrac{1}{\big\|\psi\big\|^2_2}\int_{-\infty}^{\infty}\int_{-\infty}^{\infty}{\mathcal{W}}^M_\psi \Big[f(2t-z)\Big] (b,a)\, \overline{\mathcal K^M(2t-z,a)}\, \psi_{a,b}^M(2t-z)\,da\,db\\
&=\dfrac{1}{\big\|\psi\big\|^2_2}\int_{-\infty}^{\infty}\int_{-\infty}^{\infty}\exp\left\{\dfrac{-i}{2B}\Big(-A(-2t)^2-2t(u_{0}-a)\Big) \right\}\\
&\qquad \times {\mathcal{W}}^M_\psi \Big[\exp\left(\dfrac{-2At^2}{B}\right)f(-z)\Big] (b,a)\,  \overline{\mathcal K^M(2t-z,a)}\, \psi_{a,b}^M(2t-z)\,da\,db\\
\end{align*}
\begin{align*}
&=\dfrac{1}{\big\|\psi\big\|^2_2}\int_{-\infty}^{\infty}\int_{-\infty}^{\infty}\exp\left\{\dfrac{i}{B}\Big(2At^2-2t(u_{0}-a)\Big) \right\}\\
&\qquad \times {\mathcal{W}}^M_\psi\Big[\exp\left(\dfrac{-2At^2}{B}\right)f(z)\Big] (-b-2t,-a)\,\overline{\mathcal K^M(z,a)}\, \psi_{a,b}^M(z)\,da\,db\\
&= \dfrac{1}{\big\|\psi\big\|^2_2}\int_{-\infty}^{\infty}\int_{-\infty}^{\infty}\exp\left\{\dfrac{i}{B}\Big(2At^2-2t(u_{0}-a)\Big) \right\}\\
&\qquad\times {\mathcal{W}}^M_\psi \Big[\exp\left(\dfrac{-2At^2}{B}\right)[f(z)]\Big] (-b-2t,-a)\, \psi^M_{a,b}(z)\,da\,db. \tag{5.4}
\end{align*}
By substituting (5.4) in (5.3), we obtain
\begin{align*}
&\Big[{WD}^M f\Big](t,a)\\
&= \dfrac{1}{2\sqrt {2\pi iB}} \exp\left\{\dfrac{i}{2B}\Big(-2a(Du_{0}-Dq)+D\big(a^2+u^2_{0}\big)\Big) \right\}\\
& \times \int_{-\infty}^{\infty} f(z)\,  \dfrac{1}{\big\|\psi\big\|^2_2}\int_{-\infty}^{\infty}\int_{-\infty}^{\infty}\exp\left\{\dfrac{-i}{B}\Big(2At^2-2t(u_{0}-a)\Big) \right\}\overline{{\mathcal{W}}^M_\psi \Big[\exp\left\{\dfrac{-2At^2}{B}\right\}f(z)\Big]} (-b-2t,-a)\\
&~\times \overline{ \psi^M_{a,b}(z)}\,da\,db \, \exp\left\{\dfrac{i}{2B}\Big(A(2z-2t)^2+2(2z-2t)(p-a)\Big) \right\} dz\\
&= \dfrac{1}{2\sqrt {2\pi iB}}\,\dfrac{1}{\big\|\psi\big\|^2_2} \exp\left\{\dfrac{i}{2B}\Big(-2a(Du_{0}-Dq)+D\big(a^2+u^2_{0}\big)\Big) \right\}\\
&~~\times \int_{-\infty}^{\infty}\int_{-\infty}^{\infty}\int_{-\infty}^{\infty}f(z)\,\exp\left\{\dfrac{i}{B}\Big(2Az^2+2z(u_{0}-a-2At)\Big) \right\}\, \overline{ \psi^M_{a,b}(z)}\, dz\\
&~~\times \overline{{\mathcal{W}}^M_\psi \Big[\exp\left(\dfrac{-2At^2}{B}\right)[f]\Big]} (-b-2t,-a)\,da\,db\\
&= \dfrac{1}{2\sqrt {2\pi iB}}\,\dfrac{1}{\big\|\psi\big\|^2_2}\, \exp\left\{\dfrac{i}{2B}\Big(-2a(Du_{0}-Dq)+D\big(a^2+u^2_{0}\big)\Big) \right\}\\
&~~\times \int_{-\infty}^{\infty}\int_{-\infty}^{\infty}{\mathcal{W}}^M_\psi \Big[\exp\left\{\dfrac{i}{B}\Big(2Az^2+2z(u_{0}-a-2At)\Big) \right\}\,[f]\Big](b,a)\, dz\\
&~~ \times \overline{{\mathcal{W}}^M_\psi \Big[\exp\left(\dfrac{-2At^2}{B}\right)\,[f]\Big]} (-b-2t,-a)\,da\,db.
\end{align*}
This completes the proof of the theorem. \quad \fbox

\section{ Special Affine Wave Packet Transform}

\parindent=0mm \vspace{.0in}
In pursuit of efficient representations of functions, C\'{o}rdoba and Fefferman \cite{CF} introduced the notion of wave packet systems by applying dilations, modulations and translations to the Gaussian function in the study of some classes of singular integral operators. Later on, Labate et al.\cite{Labate} adopted the same strategy to define any collections of functions which are obtained by applying the same operations to a finite family of functions in $L^2(\mathbb R)$. In fact, Gabor and wavelet systems are special cases of wave packet systems. Recent developments in this direction can be found in \cite{CR,LW,ShahA,ShahO,ShahOJ} and the references therein.

\parindent=8mm \vspace{.1in}
Following the idea of C\'{o}rdoba and Fefferman \cite{CF}, we construct a  generalized family of wave packet systems associated with the special affine Fourier transform which encompasses  all the existing wave packet systems on the Euclidean  space $\mathbb R$.

\parindent=0mm \vspace{.1in}
{\bf{Definition 6.1.}}  For a uni-modular matrix $M=(A,B,C,D,p,q)$, we define a   system of the form
\begin{align*}
{\psi}^{M}_{a,b,N}(t)=\dfrac{1}{\sqrt {2\pi iaB}}\overline{\psi} \left(\dfrac{t-b}{a}\right)\exp\left\{\dfrac{-i}{2B}\Big(At^2+2t(p-N)-2N(Dp-Bq)+D\big(N^2+p^2\big)\Big)\right\}  \tag{6.1}
\end{align*}

\parindent=0mm \vspace{.0in}
 called  special affine wave packet system  in $L^2(\mathbb R)$, where $\psi$ is a fixed function in $L^2(\mathbb R)$.

\parindent=0mm \vspace{.1in}
{\bf Definition 6.2.} Given a real unimodular matrix $M=\big(A,B,C,D,p,q\big)$, the special affine wave packet transform of any square integrable function $f$ is defined by
\begin{align*}
\Big[{WP}^M_\psi f\Big](a,b,N)=\int_{-\infty}^{\infty}f(t)\,\overline{{\psi}^{M}_{a,b,N}(t)}\,dt,\tag{6.2}
\end{align*}
where ${\psi}^{M}_{a,b,N}(t)$ is given by (6.1).

\parindent=8mm \vspace{.1in}
This definition allows us to make the following observations:

\parindent=0mm \vspace{.1in}
(i). If we consider only the modulations and translations of a single function $\psi$ in (6.1), then  we obtain the windowed special affine Fourier transform of the form
\begin{align*}
\Big[{WP}^M_\psi f\Big](b,N)=\int_{-\infty}^{\infty}f(t)\,\overline{{\psi}^{M}_{b,N}(t)}\,dt,\tag{6.3}
\end{align*}
where
\begin{align*}
\psi^M_{b,N}(t)= \dfrac{\overline{\psi}(t-b)}{i\sqrt {2\pi i B}}\exp\left\{\dfrac{-i}{2B}\Big[At^2+2t(p-N)-2N(Dp-Bq)+D(N^2+p^2)\Big] \right\}.\tag{6.4}
\end{align*}

\parindent=0mm \vspace{.0in}
(ii). Similarly, if we consider the dilations and translations of a single function $\psi$ in (6.1), then  the special affine wave packet transform (6.2) boils down to the continuous special affine wavelet transform as defined in (2.4).

\parindent=0mm \vspace{.1in}
(iii). For the matrix $M=(\cos\theta,\sin \theta,-\sin \theta,\cos\theta,0,0)$, the special affine wave packet transform (6.2) reduces to the fractional wave packet transform given in \cite{ShahOJ}.

\parindent=0mm \vspace{.1in}
(iv). For the matrix $M=(0,1,-1,0,0,0)$, the special affine wave packet transform (6.2) reduces to the classical wave packet transform given in \cite{CR}.

\parindent=8mm \vspace{.1in}
Following theorem summarises some of the fundamental properties of the special affine wave packet transform in $L^2(\mathbb R)$.

\parindent=0mm \vspace{.1in}
{\bf Theorem 6.3.}{\it If $\psi$ and $\phi$ are the wavelet functions and  $f_1$ and $f_2$ are functions which belong to $L^2(\mathbb R)$, then the  special affine wave packet transform  defined by (6.2)  satisfies the following properties:}

\parindent=0mm \vspace{.1in}
(i) {\it Linear Property:}

\parindent=0mm \vspace{.1in}
\qquad $\Big({WP}^M_\psi\Big)\big[(\alpha f_1+\beta f_2)\big](a,b,N)=\alpha\,\Big({WP}^M_\psi\Big)\big[ f_1\big](a,b,N)+\beta\Big({WP}^M_\psi\Big)\big[f_2\big](a,b,N)$.

\parindent=0mm \vspace{.1in}
(ii) {\it Anti-Linear Property:}

\parindent=0mm \vspace{.1in}
\qquad $\Big({WP}^M_{\alpha \psi+\beta \phi} \Big)\big[f\big](a,b,N)=\overline{\alpha}\,\Big({WP}^M_\psi\Big)\big[ f\big](a,b,N)+\overline{\beta}\Big({WP}^M_\phi\Big)\big[f\big](a,b,N)$.

\parindent=0mm \vspace{.1in}
(iii)  {\it Time-Shift Property:}
\begin{align*}
\Big({WP}^M_\psi\Big)\Big[f(t-k)\Big](a,b,N)=\exp \left\{\dfrac{-iACk^2}{2}+iCk(N+p)+iAkq\right\} \Big({WP}^M_\psi\Big)(a,b-k,N-Ak).
\end{align*}

\parindent=0mm \vspace{.0in}
(iv) {\it Phase-Shift Property:}
\begin{align*}
\Big({WP}^M_\psi\Big)\Big[e^{i\alpha t} f\Big](a,b,N)=\exp\left\{-\alpha (Dp-Bq)-\dfrac{iBD \alpha^2}{2}+iN\alpha B) \right\} \Big({WP}^M_\psi\Big)\big[ f\big](a,b, N-\alpha B).
\end{align*}

\parindent=0mm \vspace{.0in}
(v). {\it Joint Phase-Time Shift Property:}
\begin{align*}
&\Big({WP}^M_\psi\Big)\Big[e^{i\alpha t} f(t-k)\Big](a,b,N)\\
&\qquad\qquad=\exp \left\{\dfrac{-i(ACk^2+BD\alpha^2)}{2}+iCk(N+p-B\alpha)+iq(Ak+B\alpha)+\alpha D(N-p)\right\}\\
&\qquad\qquad\qquad\qquad\qquad\qquad\times \Big({WP}^M_\psi\Big)(a,b-k,N-Ak-\alpha B).
\end{align*}

\parindent=0mm \vspace{.0in}
(vi). {\it Moyal's Principle:}
\begin{align*}
\Big\langle {WP}^M_\psi\big[f_1\big],{WP}^M_\psi\big[f_2\big] \Big\rangle= \big\| \psi\big\|_2 ^2\, \Big\langle f_1, f_2\Big\rangle.
\end{align*}

\section{Poisson Summation Formula For The Special Affine Wavelet Transform}

One of the fundamental formulas of classical harmonic analysis is Poisson's summation formula. This formula asserts that the sum of infinite samples in the time domain of a function $f\in L^2(\mathbb R)$ is equivalent to the sum of infinite samples of $\mathscr F[f]$ in the frequency domain. The Poisson summation formula links the signal to the samples of its spectrum by
\begin{align*}
\sum_{k\in\mathbb Z}f(t+kT)=\dfrac{1}{T}\sum_{k\in\mathbb Z}\mathscr F\big[f\big]\left(\dfrac{k}{T}\right)\exp\left\{\dfrac{ikt}{T}\right\}.\tag{7.1}
\end{align*}

\parindent=0mm \vspace{.0in}
In other words, this formula expresses the fact that discretization in one domain implies periodicity in the reciprocal domain. The Poisson summation formula has been extended in various directions and has found extensive applications in various scientific fields, particularly in signal and image processing, harmonic analysis,   quantum mechanics and number theory \cite{AM,DEK,LTXW,CHH,ZH}. Keeping in view that windowed SAFT enjoys several pleasant properties due to extra  degrees of freedom, it is worthwhile to generalize the traditional Poisson summation formula from the classical Fourier domain to the windowed special affine Fourier domain. To achieve the desired objective, we shall apply a fundamental relationship between the conventional Fourier transform and the windowed special affine Fourier transform. The following is an analogue of the Poisson summation formula for the special affine wavelet transform.

\parindent=0mm \vspace{.1in}
{\bf Theorem 7.1.} {\it If ${\mathcal W}_\psi^M \big[f\big](a,b)$ is the continuous special affine wavelet transform of any $f(t)\in L^2(\mathbb R)$, then associated  Poisson summation formula  is given by}
\begin{align*}
&\dfrac{1}{\sqrt{a}}\,\sum_{k\in\mathbb Z}f(t+kT)\overline{\psi \,\left(\dfrac {t+kT-b}{a}\right)} \exp\left\{\frac{i}{2B}\Big(Ak^2T^2+2AkT t+2pkT\Big)\right\}\\
&= \dfrac{\sqrt{2i\pi B}}{T}\, \exp\left\{\dfrac{-i}{2B}\Big(At^2+2tp+D{p}^2\Big)\right\}
\\&\qquad\qquad\times\sum_{k\in\mathbb Z}\, \exp\left\{-\dfrac{iD}{2B}\left(\dfrac{Bk}{T}\right)^2+\left(\dfrac{ik}{T}\right)(Dp-Bq)+\dfrac{ikt}{T}\right\}\, {\mathcal{W}}_f^M \left(b,\,\frac{Bk}{T}\right).\tag{7.2}
\end{align*}

\parindent=0mm \vspace{.0in}

{\it Proof.} By virtue of the  continuous special affine wavelet transform (2.4), we have
\begin{align*}
&{\mathcal{W}}_\psi^M \left(\frac{Bk}{T},\,b\right)=\dfrac{1}{\sqrt{2\pi iBa}}\int_{-\infty}^{\infty}f(t)\,\overline{\psi\,\left(\dfrac {t-b}{a}\right)}\\
&\qquad\times\exp\left\{\frac{i}{2B}\Big(At^2+2t\left(p-\dfrac{Bk}{T}\right)-2\left(\dfrac{Bk}{T}\right)\left(Dp-Bq\right)
+D\left(\dfrac{Bk}{T}\right)^2+Du_o^2 \Big)\right\}dt.
\end{align*}
The above equation can also be rewritten in the following manner:
\begin{align*}
&\sqrt{2\pi iB}\,\exp\left\{-\dfrac{iD}{2B}\left(\dfrac{Bk}{T}\right)^2\right\}\, \exp\left\{\dfrac{-i}{2B}\Big(D{p}^2-2\left(\dfrac{Bk}{T}\right)\left(Dp-Bq\right)\Big)\right\}
\, {\mathcal{W}}_f^M \left(b,\,\frac{Bk}{T}\right)\\
&=\int_{-\infty}^{\infty}\exp\left\{\dfrac{i}{2B}\Big(At^2+2tp\Big)\right\} f(t)\, \overline{\psi\,\left(\dfrac {t-b}{a}\right)}\, \exp\left\{-\dfrac{ikt}{T}\right\}\,dt\\
&=\mathscr F\left[\exp\left\{\dfrac{i}{2B}\Big(At^2+2tp\Big)\right\}\, f(t)\, \dfrac{1}{\sqrt{a}}\,\overline{\psi\,\left(\dfrac {t-b}{a}\right)}\right]\left(\dfrac{k}{T}\right).
\end{align*}

\parindent=0mm \vspace{.0in}
We shall make use of the notation $h(t)=e^{\frac{i}{2B}\left(At^2+2tp\right)}\, f(t)\,\overline{\psi_{a,b}(t)}$ and applying the formula (7.1) to the function $h(t)$, we obtain
\begin{align*}
&\sum_{k\in\mathbb Z}h\left(t+kT\right)=\dfrac{1}{T}\sum_{k\in\mathbb Z}\mathscr F\left[h\right]\left(\dfrac{k}{T}\right)\, \exp\left\{\dfrac{ikt}{T}\right\}\\
&=\dfrac{\sqrt{2i\pi B}}{T}\, \exp\left(-\dfrac{iD{p}^2}{2B}\right)
\,\sum_{k\in\mathbb Z}\,\exp\left\{-\dfrac{iD}{2B}\left(\dfrac{Bk}{T}\right)^2+\left(\dfrac{ik}{T}\right)(Dp-Bq)+\dfrac{ikt}{T}\right\} {\mathcal{W}}_\psi^M \left(\frac{Bk}{T},\,b\right).\tag{7.3}
\end{align*}

\parindent=0mm \vspace{.1in}
Furthermore, we observe that
\begin{align*}
&\sum_{k\in\mathbb Z}h\left(t+kT\right)\\
&=\sum_{k\in\mathbb Z}f(t+kT)\,\dfrac{1}{\sqrt{a}} \overline{\psi\,\left(\dfrac {t+kT-b}{a}\right)} \exp\left\{\dfrac{i}{2B}\Big(A(t+kT)^2+2(t+kT)p\Big)\right\}\\
&=\dfrac{1}{\sqrt{a}}\sum_{k\in\mathbb Z}f(t+kT)\,\overline{\psi\left(\dfrac {t+kT-b}{a}\right)}\, \exp\left\{\frac{i}{2B}\Big(At^2+Ak^2T^2+2AkT t+2tp+2kT p\Big)\right\}\\
&=\dfrac{1}{\sqrt{a}}\,\exp\left\{\dfrac{i}{2B}\Big(At^2+2tp\Big)\right\}\sum_{k\in\mathbb Z}f(t+kT)\,\overline{\psi\left(\dfrac {t+kT-b}{a}\right)}\, \exp\left\{\frac{iA}{2B}\left(k^2T^2+2kT t\right)+\dfrac{ipkT}{B}\right\}.\tag{7.4}
\end{align*}

\parindent=0mm \vspace{.1in}
Using  equation (7.4) in (7.3), we get the desired result
\begin{align*}
&\dfrac{1}{\sqrt{a}}\,\sum_{k\in\mathbb Z}f(t+kT)\overline{\psi}\left(\dfrac {t+kT-b}{a}\right) \exp\left\{\frac{i}{2B}\Big(Ak^2T^2+2AkT t+2pkT\Big)\right\}\\
&= \dfrac{\sqrt{2i\pi B}}{T} \exp\left\{\dfrac{-i}{2B}\Big(At^2+2tp+D{p}^2\Big)\right\}
\\&\qquad\times\sum_{k\in\mathbb Z} \exp\left\{-\dfrac{iD}{2B}\left(\dfrac{Bk}{T}\right)^2+\left(\dfrac{ik}{T}\right)(Dp-Bq)+\dfrac{ikt}{T}\right\} {\mathcal{W}}_\psi^M \left(\frac{Bk}{T},\,b\right).
\end{align*}

\parindent=0mm \vspace{.0in}
In particular, for $t=0$, we have
\begin{align*}
&\dfrac{1}{\sqrt{a}}\sum_{k\in\mathbb Z}f(kT)\,\overline{\psi\,\left(\dfrac {kT-b}{a}\right)}\,\exp\left\{\frac{i}{2B}\Big(Ak^2T^2+2pkT\Big)\right\}\\
&= \dfrac{\sqrt{2i\pi B}}{T} \exp\left\{\dfrac{-i}{2B}\left(D{p}^2\right)\right\}\, \sum_{k\in\mathbb Z}\, \exp\left\{-\dfrac{iD}{2B}\left(\dfrac{Bk}{T}\right)^2+\left(\dfrac{ik}{T}\right)(Dp-Bq)\right\}{\mathcal{W}}_\psi^M \left(\frac{Bk}{T},\,b\right).
\end{align*}

This completes the proof of Theorem 7.1. \qquad\fbox

\end{document}